\documentclass{article}
\usepackage{graphicx} 
\usepackage{tikz}

\usepackage{amsmath}
\usepackage{amsfonts}
\usepackage{amssymb}
\usepackage{amsthm}
\usepackage{setspace}
\usepackage{hyperref}
\usepackage{graphicx}
\usepackage{stmaryrd}
\usepackage{mathtools}
\usepackage{subcaption}

\usepackage{stackengine}

\newtheorem{theorem}{Theorem}[section]
\newtheorem{proposition}[theorem]{Proposition}
\newtheorem{lemma}[theorem]{Lemma}

\theoremstyle{definition}

\title{A virtual RAAG with no finite index normal RAAG}
\author{Oli Jones \footnote{Email: jones@math.tu-berlin.de. Address: Technische Universit\"at Berlin, Straße des 17. Juni 135, 10623 Berlin, Germany.}}
\date{June 2026}

\begin{document}

\maketitle

\begin{abstract}
    In this note, we exhibit a group which has a right-angled Artin group as a finite index subgroup, but no finite index normal subgroup isomorphic to any right-angled Artin group. This answers a recent question of Vankov.
\end{abstract}

\section{Introduction}

Many natural classes of finitely generated groups $\mathcal{C}$ are closed under passing to finite index subgroups. This is true, for instance, of free groups and free abelian groups. It is easy to see that such classes enjoy a second property, that if any group $G$ is virtually $\mathcal{C}$, then $G$ has a final index normal subgroup in $\mathcal{C}$. 

Right-angled Artin groups (RAAGs) are a rich and pervasive class within geometric group theory, which interpolate between free and free abelian groups. It is not true, however, that a finite index subgroup of a RAAG is itself a RAAG. The purpose of this note is to show that RAAGs also fail the second property, answering a question of Vankov. Namely, we exhibit a group which has $A_{C_6}$ as an index 3 subgroup, but which has no finite index normal subgroup which is isomorphic to any RAAG. Note that the index of 3 is optimal, since any subgroup of index 1 or 2 is itself normal.

\section{Right-angled Artin groups}

Given a finite simplicial graph $\Gamma$, the associated RAAG is the the group generated by $V(\Gamma)$, with the relations that vertices adjacent in the graph commute. We will call the generators corresponding to the vertices standard generators of $A_\Gamma$. 

Given a subset of the generators $U \subseteq V(\Gamma)$, the subgroup of $A_\Gamma$ generated by $U$ is called a standard parabolic subgroup. It is well known that this subgroup is isomorphic to the RAAG on the subgraph spanned by $U$ \cite{van1983homotopy}. We will denote this subgroup by $A_U$. Conjugates of standard parabolic subgroups are called parabolic subgroups.

We will repeatedly exploit the centraliser structure of RAAGs. Given $g \in A_\Gamma$, the parabolic closure is the smallest (with respect to inclusion) parabolic subgroup containing $g$. Up to conjugation we can assume the parabolic closure of $g$ is a standard parabolic subgroup. The centraliser of $g$ admits a complete description in terms of this parabolic closure; we will use only the following consequence.

\begin{proposition}\label{prop:cent}\cite{servatius1989automorphisms}
    Suppose that $\Gamma$ is triangle-free and square-free. Then $g \in A_\Gamma$ has a non-abelian centraliser only if $g$ is conjugate to a power of a standard generator. Moreover, if the centraliser of $g$ is non-abelian then it contains a non-abelian free group.
\end{proposition}
\begin{proof}
    Up to conjugation, we may suppose that the parabolic closure of $g$ is a standard parabolic subgroup of $A_\Gamma$, generated by $U \subseteq A_\Gamma$.

    The centraliser of $g$ is described completely by \cite[Centralizer Theorem]{servatius1989automorphisms}. It splits as a direct product of cyclic factors, and the standard parabolic subgroup generated by $$S = \{ v \in V(\Gamma) \mid \forall u \in U, [v,u] = 1\}.$$ Evidently, the centraliser of $g$ is non-abelian if and only if $A_S$ is non-abelian. If $|U| = 1$ then $g$ is, up to conjguation, a power of a standard generator. Henceforth we restrict to the case where $|U| \geq 2$, and our aim is to show that $A_S$ is abelian. 

     Since $\Gamma_S$, the subgraph of $\Gamma$ spanned by $S$, is a join of the subgraphs $\Gamma_{S \cap U}$ and $\Gamma_{S \setminus U}$, $A_S$ splits as a direct product $A_S = A_{S \cap U} \times A_{S \setminus U}$. Now, notice that $A_{U \cap S}$ is abelian, since it is contained in the centre of $A_U$ by definition of $S$. So the centraliser of $g$ is non-abelian if and only if $A_{S \setminus U}$ is.

     Now take $u_1, u_2 \in U$ arbitrary distinct vertices (we may do this since we have that $|U| \geq 2$). There are two cases to consider. If $u_1$ and $u_2$ span an edge, then $S \setminus U$ must be empty, as any vertex in $S \setminus U$ would form a triangle with $u_1$ and $u_2$ by definition of $S$. Otherwise, we have $|S \setminus U| \leq 1$, since if $S \setminus U$ contained 2 distinct vertices, they would span a square with $u_1$ and $u_2$. So $A_{S \setminus U}$ is cyclic or trivial and, in particular, abelian. 

    The moreover part of the statement follows because every non-abelian subgroup of a RAAG contains a non-abelian free group \cite{baudisch1981subgroups}.
\end{proof}

\section{The example}

The example will be a group $$G = \langle t \mid t^3 \rangle \rtimes A_{C_6},$$ where $$A_{C_6} = \langle a_1 \dots a_6 \mid [a_i, a_{i+1}]\rangle,$$ is the hexagon RAAG, and by convention we take $a_7 = a_1$. Each standard generator of $A_{C_6}$ will act on $\mathbb{Z}/3\mathbb{Z}$ by the sole non-trivial automorphism. That is, for each $i \in \{1 \dots 6\}$, we impose that $a_ita_i^{-1} = t^2$.

Evidently, $A_{C_6}$ is an index $3$ subgroup of $G$. The rest of this section is devoted to a proof of the following.

\begin{theorem}\label{thm:main}
    The group $G$ does not have a finite index normal subgroup isomorphic to a RAAG.
\end{theorem}

In fact, the crux of the argument is the following theorem.

\begin{theorem}\label{thm:ac6}
    Suppose $N \leq A_{C_6}$ is a finite index normal subgroup isomorphic to a RAAG. Then, for some $i$, $a_i \in N$.
\end{theorem}
\begin{proof}
    Suppose that $N \leq A_{C_6}$ is a finite index normal subgroup isomorphic to a RAAG. We fix an identification $N = A_\Gamma$, so we may talk about (for example) standard generators of $N$ to mean the elements corresponding to vertices of $\Gamma$. Note that $\Gamma$ must be connected, as if not then $A_\Gamma$ is infinitely-ended, but $A_{C_6}$ is one-ended. One-endedness is a quasi-isometry invariant and in particular passes to finite index subgroups. The graph $\Gamma$ must have at least 3 vertices, since $A_{C_6}$ is not virtually abelian.

    The group $A_{C_6}$ has cohomological dimension 2. It follows that $N$ must also have cohomological dimesnion 2. In particular, $\Gamma$ is triangle-free. Moreover, $A_{C_6}$ is square-free, and since $A_\Gamma$ is commensurable to $A_{C_6}$ it must also be square-free by \cite[Corollary 1.6]{genevois2026polynomialhyperbolicityproductsfree}. As such, we freely apply Proposition \ref{prop:cent} to $A_{C_6}$ and $N = A_\Gamma$.

    The next two lemmas will restrict which elements of $A_{C_6}$ may appear as standard generators of $N$. We do this by showing that such elements' centralisers must contain a free subgroup. Note that $g \in N$ has a free subgroup in its centraliser in $N$ if and only if it has a free subgroup in its centraliser in $A_{C_6}$, so we may freely pass this property between the two groups.

    \begin{lemma}\label{lem:nonab}
        Every standard generator of $N$ has a free subgroup in its centraliser in $N$.
    \end{lemma}
    \begin{proof}
        Since $\Gamma$ is triangle-free, it follows that a standard generator $v \in V(\Gamma)$ may only have abelian centraliser if it is a leaf of $\Gamma$. 

        Now, however, recall that a RAAG splits over $\mathbb{Z}$ if and only if the defining graph has a separating vertex \cite{clay2014does}, so $A_{C_6}$ does not split over $\mathbb{Z}$. However $A_{C_6}$ and $A_\Gamma$ are, in particular, quasi-isometric. Splitting over $\mathbb{Z}$ is a quasi-isometry invariant amongst finitely presented torsion-free groups \cite{papasoglu2005quasi}, so $A_\Gamma$ does not split over $\mathbb{Z}$, and hence $\Gamma$ has no separating vertex. Since $\Gamma$ has at least 3 vertices, this means it cannot have leaves. The result follows. 
    \end{proof}

    \begin{lemma}\label{lem:formofSG}
        If $g \in N \leq A_{C_6}$ is a standard generator of $N$, then $g$ is conjugate in $A_{C_6}$ to $a_i^n$, where $n \in \mathbb{Z}^*$.
    \end{lemma}
    \begin{proof}
        By Lemma \ref{lem:nonab}, $g$ has a free subgroup in its centraliser in $N$. Clearly it follows that $g$ has a free subgroup in its centraliser in $A_{C_6}$, and the result follows by Proposition \ref{prop:cent} applied to $A_{C_6}$.
    \end{proof}

    We are now ready to prove the theorem. For each $i \in \{1 \dots 6\}$, fix $n_i \in \mathbb{N}^*$ to be minimal such that $a^{n_i} \in N$. Such $n_i$ must exist, since $N$ has finite index in $A_{C_6}$. Clearly, $a_i^k \in N$ if and only if $k$ is a multiple of $n_i$. Since $N$ is normal, we in particular have that $ga_1^{n_1}g^{-1} \in N$ for each $g \in G$.

    For each $g \in A_{C_6}$, note that the centraliser of $ga_1^{n_1}g^{-1}$ in $A_{C_6}$ contains a free subgroup. It follows that the centraliser of $ga_1^{n_1}g^{-1}$ in $N$ contains a free subgroup and, hence, by Proposition \ref{prop:cent} applied to $N$, $ga_1^{n_1}g^{-1}$ must be conjugate in $N$ to a power of a standard generator of $N$. In fact it must be conjugate exactly to a standard generator or its inverse, since $ga_1^{n_1}g^{-1}$ is root-free in $N$ by the minimality of $n_1$. Since $N$ has only finitely many standard generators, it follows that $$\{ga_1^{n_1}g^{-1} \mid g \in A_{C_6}\},$$ meets only finitely many conjugacy classes in $N$. 

    We now fix, for each $m \in \mathbb{N}^*$, $$w_m = (a_5a_3)^m,$$ to be the alternating word of length $2m$. By the discussion above, we may find $m, k \neq 0$ such that $w_ma_1^{n_1}w_m^{-1}$ and $w_{m+k}a_1^{n_1}w_{m+k}^{-1}$ are conjugate in $N$. Take $g \in N$ to be such that $$gw_ma_1^{n_1}w_m^{-1}g^{-1} = w_{m+k}a_1^{n_1}w_{m+k}^{-1},$$ or equivalently such that $w_{m+k}^{-1}gw_m$ centralises $a_1^{n_1}$ in $A_{C_6}$. It follows by \cite{servatius1989automorphisms} that the centraliser of $a_1^{n_1}$ in $A_{C_6}$ is given by $\langle a_1, a_2, a_6 \rangle$. In particular, if $\phi: A_{C_6} \rightarrow \langle a_3, a_5 \rangle$ is the retraction sending each other standard generator to $1$, we can apply $\phi$ to conclude that $$w_{m+k}^{-1}\phi(g)w_m = 1,$$ and hence that $\phi(g) = w_k$.

    We now use that we have $g \in N$, and so $g$ can be written as a product of standard generators of $N$. By Lemma \ref{lem:formofSG}, therefore, we have that $$g \in \langle \langle \{ a_i^{n_i} \mid i \in \{1 \dots 6\}\} \rangle \rangle,$$ and hence that, $$w_k = \phi(g) \in \langle \langle a_3^{n_3}, a_5^{n_5} \rangle \rangle.$$ However the right hand side is the kernel of the map from $\langle a_3, a_5 \rangle$ to the free product $$\langle a_3, a_5 \mid a_3^{n_3}, a_5^{n_5} \rangle,$$ and by the normal form for this free product we see that we must have $n_3 = 1$ or $n_5 = 1$.

\end{proof}

We are now ready to prove Theorem \ref{thm:main}. We will make use of the obvious retraction $\psi: G \rightarrow A_{C_6}$ arising from the semi-direct product structure, whose kernel is the finite subgroup $K = \langle t \mid t^3 \rangle$.

\begin{proof}[Proof of Theorem \ref{thm:main}]
    Suppose for contradiction that $N \leq G$ is a finite index normal subgroup isomorphic to a RAAG. Then it must be that $K \cap N = \{1\}$, since $K$ is finite and $N$ is torsion-free. Hence we have that $\psi(N) \cong N$ is isomorphic to a RAAG. Clearly $\psi(N) \leq A_{C_6}$ is a finite index subgroup, since it contains $N \cap A_{C_6}$ which is finite index in $A_{C_6}$. Moreover, if $g \in A_{C_6}$, then $$g\psi(N)g^{-1} = \psi(gNg^{-1}) = \psi(N),$$ so $\psi(N)$ is normal in $A_{C_6}$.

    Now, by Theorem \ref{thm:ac6}, take $a_i$ a standard generator of $A_{C_6}$, such that $a_i \in \psi(N)$. Then can find $g \in \psi^{-1}(a_i) \cap N$, and notice that $g = a_it^k$ for $k \in \{0,1,2\}$. Since $N$ is a normal subgroup, we have $$N \ni g^{-1}(tgt^{-1}) = t^{-k}a_i^{-1}ta_it^{k-1} = t^{-k}t^2t^{k-1}= t,$$ but this contradicts that $N$ is torsion-free.
\end{proof}

\section*{Acknowledgements}

I am grateful to Alex Levine for bringing this question back to my attention, and for spotting a mistake in the first version of the proof. I am grateful to Vlad Vankov for asking this question and encouraging me to write up the result. Finally, thanks to Laura Ciobanu, Sobhi Massalha and Daniel Turaev for their useful comments.

\bibliographystyle{amsalpha}
\bibliography{bibliography}

\end{document}